\documentclass[11pt]{amsart}
\usepackage[paperwidth=17cm, paperheight=22.5cm, bottom=2.5cm, right=2.5cm]{geometry}
\usepackage{amssymb,amsmath,amsthm} 
\usepackage[english]{babel}
\usepackage[utf8x]{inputenc}
\usepackage[T1]{fontenc}
\usepackage{enumerate}
\usepackage{mathtools}
\usepackage{graphicx}
\usepackage{subfigure} 
\usepackage{float}
\graphicspath{{Img/}} 
\usepackage[export]{adjustbox}
\usepackage{tikz}
\usetikzlibrary{matrix}
\usepackage{multirow}

\newtheorem{teo}{Theorem}[section]
\newtheorem{propo}[teo]{Proposition}
\newtheorem{cor}[teo]{Corollary}
\newtheorem{lem}[teo]{Lemma}
\newtheorem{que}{Question}
\newtheorem*{que*}{Question}

\theoremstyle{definition}
\newtheorem{dfn}[teo]{Definition}

\theoremstyle{remark}

\newcommand{\conc}{^\smallfrown}	
\newcommand{\rest}{\upharpoonright} 
\newcommand{\forces}{\Vdash}		

\newcommand{\QED}{\hfill\ensuremath{\square}}

\def\cA{{\mathcal{A}}} \def\cB{{\mathcal{B}}}    \def\cF{{\mathcal{F}}}   \def\cI{{\mathcal{I}}} \def\cJ{{\mathcal{J}}}      \def\cP{{\mathcal{P}}}        \def\cX{{\mathcal{X}}}

\title{Submaximal spaces and cardinal invariants}
\author{C\'esar Corral}
\address{Centro de Ciencias  Matem\'aticas, Universidad Nacional Aut\'onoma de M\'exico,
\' Campus Morelia, 58089, Morelia, Michoac\'an, M\'exico.}
\email{cicorral@matmor.unam.mx}

\date{}

\thanks{%
The author gratefully acknowledges support from CONACyT scholarship 742627.}

\keywords{}
\subjclass[2020]{54A25,54A35, 03E17,03E35}
\begin{document}

\begin{abstract}
We give a combinatorial characterization of countable submaximal subspaces of $2^\kappa$. Using a parametrized version of Mathias forcing, we prove that there exists a countable submaximal subspace of $2^{\omega_1}$ whilst $\mathfrak{c}=\omega_2$. Combining this with previous results, we construct a disjointly tight countable irresolvable space of weight $<\mathfrak{c}$, answering a question of Bella and Hru\v{s}\'{a}k.
\end{abstract}

\maketitle

\section{Introduction}

All spaces considered are Hausdorff and crowded, i.e., with no isolated points. A space $X$ is \emph{irresolvable} if for every dense $D\subseteq X$ the complement $X\setminus D$ is not dense \cite{hewitt1943problem}. A space $X$ is \emph{submaximal} iff every dense subspace $D$ is also open \cite{hewitt1943problem}, hence, submaximal spaces are trivially irresolvable. 


Given a set $\cX$, we denote the set of all finite partial functions $p:\cX\to 2$ by $FF(\cX)$.
Recall that a family $\cA\subseteq\cP(\omega)$ is \emph{independent} if for every $h\in FF(\cA)$, the set $\cA^h:=\bigcap_{A\in dom(h)}A^{h(A)}$, is infinite. Here, $X^0=X$ and $X^1=\omega\setminus X$ for every $X\subseteq\omega$. Each time we write $\cA^h$ for an independent family $\cA$, we will be referring to this kind of sets.
We denote by $\mathfrak{i}$ the minimal size of a maximal independent family. 
Our notation is mainly standar and follows \cite{engelking1989general} and \cite{kunen2014set}. In particular, given an ideal $\cI\subseteq\cP(\omega)$, its dual filter is $\cI^*=\{\omega\setminus X:X\in\cI\}$ and similarly we define by $\cF^*$ the dual ideal of a filter $\cF\subseteq\cP(\omega)$. The set of positive sets with respect to an ideal $\cI$ is $\cI^+=\cP(\omega)\setminus\cI$ and for a filter $\cF$, we write $\cF^+:=(\cF^*)^+$.

Given any forcing notion $\mathbb{P}$ and a $\mathbb{P}$-generic filter $G$ over the universe $V$, we denote by $V[G]$ its extension by $G$ and $V^\mathbb{P}$ denotes any arbitrary extension of $V$ after forcing with $\mathbb{P}$, without specify the generic filter. If $\mathbb{P}=\langle\mathbb{P}_\alpha,\dot{\mathbb{Q}}_\beta:\alpha\leq\kappa,\beta<\kappa\rangle$ is a finite support iteration of forcing notions, we will assume that its elements $p\in\mathbb{P}$ are finite functions with domain a subset of $\kappa$ and such that $(p\rest\alpha)\forces_\alpha``p(\alpha)\in\dot{\mathbb{Q}}_\alpha$'' for every $\alpha\in dom(p)$. We will also denote by $supp(p)$ the domain of $p$.

\section{Cardinal invariants associated to submaximal spaces}

The study of resolvable spaces began with the work of Hewitt in 1943 \cite{hewitt1943problem}. There, a space is called \emph{resolvable} if it contains two disjoint dense subsets and it is irresolvable otherwise. Resolvable and irresolvable spaces have been widely studied and it is known that some classes of spaces are resolvable, e.g., locally compact Hausdorff spaces, metrizable spaces and countably compact regular spaces. It is intriguing that it is not known if pseudocompact Tychonoff spaces are resolvable in {\sf ZFC} (they are under $V=L$!). 
On the other hand, several examples of irresolvable spaces have been constructed, but the properties of these examples are strongly tied to the behaviour of cardinal invariants. For more on resolvable an irresolvable spaces see \cite{comfort1996resolvability}.

In \cite{scheepers2007topological}, the cardinal $\mathfrak{irr}$ is defined as the minimal $\pi$-weight of a countable regular irresolvable space. 
In \cite{cancino2014countable} it is proved that the cardinal $\mathfrak{irr}$ is also equal to the minimum size of a base for a countable regular irresolvable space, i.e., we can exchange $\pi$-weight and weight in the definition of $\mathfrak{irr}$.
We can define the analogue cardinal invariant for submaximal spaces:

\begin{dfn}
$$\mathfrak{sbm}=\min\big\{|\cA|:\cA\textnormal{ is a $\pi$-base for a $T_3$ submaximal topology on $\omega$}\big\}.$$
\end{dfn}

It is known that one can construct a countable dense irresolvable space in $2^{|\cA|}$ whenever $\cA$ is maximal independent, thus $\mathfrak{irr}\leq\mathfrak{i}$ \cite{franco2003topological}. A lower bound for $\mathfrak{irr}$ is also known.

\begin{propo}\cite{scheepers1999combinatorics}
$\mathfrak{r}\leq\mathfrak{irr}$.
\end{propo}

\begin{proof}
Let $\kappa<\mathfrak{r}$ and $\{U_\alpha:\alpha<\kappa\}$ be a family of open sets of a countable regular irresolvable space $X=(\omega,\tau)$. Thus, $\{U_\alpha:\alpha<\kappa\}$ is a family of less than $\mathfrak{r}$ infinite subsets of $\omega$. Let $B\subseteq\omega$ be a set that splits every $U_\alpha$. If $\{U_\alpha:\alpha<\kappa\}$ were a $\pi$-base, then $B$ would be a dense and co-dense subset of $X$. Therefore, $\{U_\alpha:\alpha<\kappa\}$ is not a $\pi$-base.
\end{proof}

It should be evident that dense countable subspaces of cantor cubes $2^\kappa$ play a central role in the study of irresolvable and submaximal spaces. As we have seen, a maximal independent family $\cA$ gives a countable dense irresolvable subspace of $2^{|\cA|}$. This gives a topological characterization of the cardinal invariant $\mathfrak{i}$. Isolating the combinatorial properties of countable dense subspaces of $2^\kappa$ which are submaximal, a natural strengthening of $\mathfrak{i}$ arises, but before we define it, let us introduce some notation. 

For an independent family $\cA$, we denote by $\cB[\cA]$ the boolean algebra generated by $\cA$. In particular, $\cA^h\in\cB[\cA]$ for every $h\in FF(\cA)$. We say that $\cA$ separates points, if for every $n,m\in\omega$, there exists $A\in\cA$ such that either $n\in A$ and $m\notin A$ or $n\notin A$ and $m\in A$.

\begin{dfn}
Let $\cA\subseteq\cP(\omega)$ be an independent family. We say that $X\subseteq\omega$ is \emph{$\cA$-dense}, if $X\cap \cA^h$ is infinite for every $h\in FF(\cA)$.\\
We will say that $\cA$ is an \emph{anchored independent family} if for every $\cA$-dense set $X$ and every $n\in X$, there exist $\cA^h\in\cB[\cA]$ such that 
$$n\in\cA^h\subseteq^* X.$$
Let $\mathfrak{i}_{a}$ be the minimum size of a anchored independent family which separates points.
\end{dfn}

It follows that $\mathfrak{i}\leq\mathfrak{i}_{a}\leq\mathfrak{c}$ from the following results.

\begin{lem}
Every anchored independent family is maximal.
\end{lem}

\begin{proof}
Let $\cA$ be an anchored independent family and let $X\subseteq\omega$. If there is $\cA^h\in\cB[\cA]$ such that $A^h\cap X$ is finite, we are done. Otherwise, $X$ is $\cA$-dense and there is an $\cA^h\in\cB[\cA]$ such that $\cA^h\subseteq^* X$. Hence $(\omega\setminus X)\cap\cA^h$ is finite and $\cA$ is maximal.
\end{proof}

\begin{propo}\label{submaximalequivalence}
The cardinal $\mathfrak{i}_{a}$ is the minimum $\kappa$ for which there is a countable, dense subspace of $2^\kappa$ which is submaximal.
\end{propo}

\begin{proof}
Let $\cA=\{A_\alpha:\alpha<\kappa\}$ be an anchored independent family which separates points and such that $\kappa=\mathfrak{i}_{a}$. For every $n\in\omega$ define $x_n\in2^\kappa$ by 
$$x_n(\alpha)=\begin{cases} 0 & \mbox{if }n\in A_\alpha\\
1 & \mbox{otherwise}\\
\end{cases}$$
Hence, $X=\{x_n:n\in\omega\}$ is a dense countable subspace of $2^\kappa$, since $\cA$ is independent. Notice that $x_n\neq x_m$ whenever $n\neq m$, since $\cA$ separates points. Note also that $X$ is Hausdorff and $\{U_\alpha,X\setminus U_\alpha:\alpha<\kappa\}$, where $U_\alpha=\{x_n:x_n(\alpha)=0\}$, is a subbase consisting of clopen sets. Hence $X$ is also zero dimensional and in consequence completely regular.
Let $Y\subseteq X$ be a dense subset and $x\in Y$. Thus $Y'=\{n\in\omega:x_n\in Y\}$ is $\cA$-dense, and if $x=x_n$, we have that $n\in Y'$. By the definition of an anchored family, there is an $\cA^h\in\cB[\cA]$ such that $n\in\cA^h\subseteq^*Y$. Then, 
$$x_n\in V=\Big(\bigcap\{U_\alpha:h(A_\alpha)=0\}\Big)\cap\Big(\bigcap\{(X\setminus U_\alpha):h(A_\alpha)=1\}\Big)\subseteq^*Y.$$
Thus we can find an open set $U$ such that $x_n\in U\subseteq V\cap Y$, which proves that $Y$ is open.

On the other hand, if $X=\{x_n:n\in\omega\}\subseteq 2^\kappa$ is a dense subspace, we can define $A_\alpha=\{n\in\omega:x_n(\alpha)=0\}$, and then $\cA=\{A_\alpha:\alpha<\kappa\}$ is independent by density of $X$. Let $n\neq m\in\omega$ and find $\alpha<\kappa$ such that $x_n(\alpha)\neq x_m(\alpha)$, then the set $A_\alpha$ separates $n$ and $m$. Moreover, if $Y\subseteq \omega$ is $\cA$-dense, $Y'=\{x_n:n\in Y\}$ is dense in $X$. Let $n\in\omega$ and let $U$ be a basic open set such that $x_n\in U\subseteq Y'$. Assume $p;\kappa\to2$ is a finite function such that
$$U=\{x\in X:\forall \alpha\in dom(p)(x(\alpha)=p(\alpha))\}.$$
Define $h\in FF(\cA)$ such that $dom(h)=\{A_\alpha:\alpha\in dom(p)\}$ and for every $\alpha\in dom(h)$, $h(A_\alpha)=p(\alpha)$.
It is easy to see that
$$n\in\cA^h\subseteq Y.$$
Therefore $\cA$ is an anchored independent family.
\end{proof}

Finally, with the next result we set an upper bound for $\mathfrak{i}_a$.

\begin{teo}\cite{juhasz2006d}
There exists a countable, dense, submaximal subspace of $2^{\mathfrak{c}}$. 
\end{teo}

Figure \ref{diagram} shows the relationship between the cardinal invariants associated to irresolvable and submaximal spaces, where an arrow from $a$ to $b$, means $a\leq b$.

{\begin{figure}[t]\label{diagram}
    \centering
    \begin{tikzpicture}
      \matrix (m) [matrix of math nodes, row sep=2em,
        column sep=2em]{
         & \mathfrak{i} & \mathfrak{i}_a & \mathfrak{c}\\
         \mathfrak{r} & \mathfrak{irr} & \mathfrak{sbm} & \\
};
    \path[-stealth]
        (m-2-1) edge (m-2-2)
        (m-2-2) edge (m-1-2) edge (m-2-3) 
        (m-2-3) edge (m-1-3)
        (m-1-2) edge (m-1-3)
        (m-1-3) edge (m-1-4);
    \end{tikzpicture}
    \caption{Cardinal invariants associated to irresolvable and submaximal spaces.}
\end{figure}
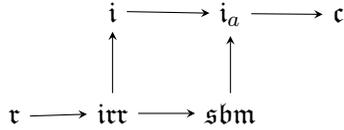}

In \cite{bella1998tightness}, the authors studied the relationship between countable tightness and resolvability. For this purpose, they introduced several variants of tightness. A space $X$ is \emph{disjointly tight} if for every $x\in X$ which is an accumulation point of a set $A$, there are disjoint $A_0,A_1\subseteq A$, such that $x\in\overline{A_0}\cap\overline{A_1}$. There, the authors prove that there exists a countable Hausdorff irresolvable spaces with disjoint tightness, and with the help of {\sf CH}, the result is improved to get a regular space. The existence of a regular space with this properties in {\sf ZFC} is left open.

In \cite{bella2020disjointly}, the authors proved that there exists a disjointly tight regular irresolvable space in {\sf ZFC}, by proving the next result:

\begin{teo}\cite{bella2020disjointly}\label{disjointly}
Let $\kappa\leq\mathfrak{c}$. Then every countable, dense, submaximal subspace of $2^\kappa$ is disjointly tight. 
\end{teo}

Thus, the existence disjointly tight irresolvable space, follows from the existence of a countable, dense, submaximal subspace of $2^\mathfrak{c}$. It has weight $\mathfrak{c}$, clearly, in consequence, the next question arose:

\begin{que*}\cite{bella2020disjointly}
Is it consistent that there is a countable, irresolvable space of weight less than $\mathfrak{c}$ that is disjointly tight?
\end{que*}

We will answer this question in the affirmative by proving that, consistently, $\mathfrak{i}_a<\mathfrak{c}$ and applying the same arguments as above.

\section{Submaximal spaces of small weight--Small anchored independent families-- \textsf{Con}($\mathfrak{i}_a<\mathfrak{c}$)}

We will make use of a parametrized version of Mathias forcing with the aim of proving that $\mathfrak{i}_a$ can be small compared with the continuum. The following lemma appears in \cite{halbeisen2012combinatorial} and \cite{fischer2019ideals}. We will give a proof in order to emphasize one property of the ideal defined.

\begin{lem}
Given an independent family $\cA\subseteq\cP(\omega)$, there exists an ideal $\cI_\cA$ such that:
\begin{enumerate}
    \item $\cI_\cA\cap\cB[\cA]=\emptyset$
    \item $\forall X\subseteq\omega\ \exists\cA^h\in\cB[\cA]\ \big((\cA^h\setminus X\in\cI_\cA)\lor(\cA^h\cap X\in \cI_\cA)\big)$
\end{enumerate}
\end{lem}

\begin{proof}
Enumerate $[\omega]^\omega=\{X_\alpha:\alpha<\mathfrak{c}\}$. Define $\cI_\alpha$ for every $\alpha<\mathfrak{c}$ as follows:
Define $\cI_0=[\omega]^{<\omega}$ and for a limit ordinal $\alpha\in\mathfrak{c}$ let $\cI_\alpha=\bigcup_{\beta<\alpha}\cI_\beta$.
Now let $\alpha=\xi+1$. If there exist an $\cA^h\in\cB[\cA]$ and $I\in\cI_\xi$ such that $\cA^h\subseteq I\cup X_\xi$, define $\cI_\alpha=\cI_\xi$. Otherwise, let $\cI_\alpha=\langle\cI_\xi\cup\{X_\xi\}\rangle$. Finally let $\cI_\cA=\bigcup_{\alpha<\mathfrak{c}}\cI_\alpha$.

To see that this works, let $X\in[\omega]^\omega$. There exists $\alpha\in\mathfrak{c}$ such that $X=X_\alpha$. At step $\alpha$, if there exist $\cA^h\in\cB[\cA]$ and $I\in\cI_\alpha$ such that $\cA^h\subseteq I\cup X_\alpha$, then $\cA^h\setminus X_\alpha\subseteq I\in\cI_\alpha\subseteq\cI_\cA$. Otherwise, $\cA^h\cap X_\alpha\subseteq X_\alpha\in\cI_{\alpha+1}\subseteq\cI_\cA$.
\end{proof}

Following \cite{fischer2019ideals}, we will call an ideal $\cI_\cA$ satisfying the previous properties, a \emph{diagonalization ideal for} $\cA$.

Notice that the ideal defined in the previous lemma is not unique.
Moreover, if $\cI'$ is an ideal and $\cI'\cap\cB[\cA]=\emptyset$, we can start the previous proof with $\cI_0=\cI'$. 

\begin{cor}\label{CoroDiagonalization ideal}
Let $\cA$ be an independent family and let $\cI$ be an ideal such that $\cI\cap\cB[\cA]=\emptyset$. There exists a diagonalization ideal $\cI_\cA$ for $\cA$ such that $\cI\subseteq\cI_\cA$. 
\end{cor}

\begin{dfn}
Let $\cA$ be an independent family and $\cI_\cA$ a diagonalization ideal for $\cA$. We define the poset $\mathbb{M}(\cI_\cA)$ as the set of all pairs $(s,E)$ such that $s\in[\omega]^{<\omega}$ and $E\in\cI_\cA$. The order is defined by $(t,F)\leq(s,E)$ if $s\subseteq t$, $E\subseteq F$ and $(t\setminus s)\cap E=\emptyset$.
\end{dfn}

Since $\mathbb{M}(\cI_\cA)$ is $\sigma$-centered, it is ccc and preserves cardinals. 
If $G$ is a $\mathbb{M}(\cI_\cA)$-generic filter, then the generic real $r_G=\bigcup\{s:\exists E\in\cI_\cA((s,E)\in G)\}\in[\omega]^\omega$ is a pseudointersection of the dual filter $\cI_\cA^*$. 
It also satisfies that $\cA\cup\{r_G\}$ is an independent family in $V[G]$.
Moreover, for every $X\in([\omega]^\omega\cap V)\setminus \cA$, the family $\cA\cup\{X,r_G\}$ is not independent and we can adjoin a maximal independent family with an iteration of posets of the form $\mathbb{M}(\cI_\cA)$. To the best knowledge of the author,the first appearance of the following theorem in the literature showed up in \cite{balcar2004combinatorics}. More detailed proofs appear in \cite{fischer2019ideals} and \cite{halbeisen2012combinatorial}.

\begin{teo}{(see \cite{balcar2004combinatorics},\cite{fischer2019ideals} and \cite{halbeisen2012combinatorial})}\label{independence}
Let $\langle\mathbb{P}_\alpha,\dot{\mathbb{Q}}_\beta:\alpha\leq\delta,\beta<\delta\rangle$ be an iteration of uncountable cofinality of posets of the form $\mathbb{M}(\cI_{\cA_\alpha})$, where $\cA_0$ is an independent family in $V$ and $\cA_\alpha$ is the independent family defined by $\cA_0\cup\{r_\beta:\beta<\alpha\}$ with $r_\beta$ the generic real adjoined at stage $\beta$. Then if $G$ is a $\mathbb{P}_\delta$-generic filter, the family $\cA_0\cup\{r_\beta:\beta<\delta\}$ is a maximal independent family in $V[G]$.\QED
\end{teo}

In order to modify the previous proof and make it work for $\mathfrak{i}_a$, we need to \emph{anchored} an element $A$ of the independent family $\cA$ to each $n\in X$  (this means, $x\in A\subseteq^* X$) whenever $X$ is $\cA$-dense. Our strategy will be to generically add elements to the independent family, that are almost contained in $X$, if possible, and otherwise, to kill the density of $X$ with respect to the independent family constructed so far. For this purpose, we need to predict those elements of the ground model that will be dense with respect to the independent family defined in the extension. The next lemma provides us of a internal way of identifying these sets.

\begin{lem}\label{A-dense}
Let $\cA$ be an independent family, $\cI_\cA$ a diagonalization ideal for $\cA$ and $X\in V\cap[\omega]^\omega$. Let $\sigma_G$ be the canonical name for the generic real added by $\mathbb{M}(\cI_\cA)$. Then $1_{\mathbb{M}(\cI_\cA)}\forces``\Check{X} \textnormal{ is }(\cA\cup\{\sigma_G\})\textnormal{-dense''}$ if and only if $(X\cap \cA^h\in\cI_\cA^+)$ for all $\cA^h\in\cB[\cA]$.
\end{lem}

\begin{proof}
$(\Rightarrow)$ Suppose that there is $\cA^h\in\cB[\cA]$ such that $X\cap \cA^h\in \cI_\cA$. Let $(s,F)\in\mathbb{M}(\cI_\cA)$, hence $(s,F\cup (X\cap \cA^h))\leq (s,F)$ and 
$$(s,F\cup(X\cap \cA^h))\forces``\sigma_G\cap(X\cap \cA^h)\subseteq s\textnormal{''.}$$
By density, $1_{\mathbb{M}(\cI_\cA)}\forces``\Check{X} \textnormal{ is not }(\cA\cup\{\sigma_g\})\textnormal{-dense''}$, a contradiction.

$(\Leftarrow)$ Suppose otherwise. Then there are $(s,F)\in\mathbb{M}(\cI_\cA)$, $\cA^h\in\cB[\cA]$, $i\in2$ and $n\in\omega$ such that
$$(s,F)\forces``|\check{X}\cap(\cA^h\cap \sigma_g^i)|<n\textnormal{''.}$$
Since $Y=(X\cap\cA^h)\setminus F$ is infinite we can find $m\in\omega$ such that $Y\cap(\max(s),m]>2n$. Take $a\in [Y\cap(\max(s),m]]^n$ and let $b=Y\cap(\max(s),m]\setminus a$. Then for $i\in2$ we have that $(s\cup a,F\cup b)\leq(s,F)$ and
$$(s\cup a,F\cup b)\forces``|X\cap\cA^h\cap\sigma^i_G|\geq n\textnormal{''.}$$
contradicting the choice of $(s,F)$.
\end{proof}

\begin{dfn}
Let $\cA$ be an independent family and $\cI_\cA$ a diagonalization ideal for $\cA$.
\begin{itemize}
    
    \item A set $X\subseteq\omega$ is $\cI_\cA$\emph{-big} if $X\cap \cA^h\in\cI_\cA^+$ for all $\cA^h\in\cB[\cA]$.
    
    \item $\cI_\cA$ is a \emph{complete diagonalization ideal for} $\cA$, if $(\omega\setminus X)\in\cI_\cA$ for every $X\subseteq\omega$ that is $\cI_\cA$-big.
\end{itemize}
\end{dfn}

\begin{lem}\label{completeideal}
Let $\cA$ be an independent family and $\cI_\cA$ a diagonalization ideal for $\cA$. There exists a complete diagonalization ideal $\cJ_\cA$ for $\cA$ such that $\cI_\cA\subseteq\cJ_\cA$.
\end{lem}

\begin{proof}
Enumerate $\{X\subseteq\omega:X \textnormal{ is }\cI_\cA\textnormal{-big}\}$ as $\{X_\alpha:\alpha<\mathfrak{c}\}$. We will recursively define $\cJ_\alpha$ for $\alpha<\mathfrak{c}$. If $\alpha=0$, simply define $\cJ_0=\cI_\cA$. For $\alpha$ a limit ordinal define $\cJ_\alpha=\bigcup_{\beta<\alpha}\cJ_\beta$. Assume $\alpha=\beta+1$. If there exist $\cA^h\in\cB[\cA]$ and $u\in\cJ_\beta$ such that $\cA^h\subseteq u\cup(\omega\setminus X_\beta)$ define $\cJ_\alpha=\cJ_\beta$. Otherwise, let $\cJ_\alpha=\langle\cJ_\beta\cup\{\omega\setminus X\}\rangle$. Finally define $\cJ_\cA=\bigcup_{\alpha<\mathfrak{c}}\cJ_\alpha$.

It is clear that $\cI_\alpha\subseteq\cJ_\alpha$. In particular, for all $X\in[\omega]^\omega$, there exists $\cA^h\in\cB[\cA]$ such that either, $\cA^h\setminus X\in\cI_\cA\subseteq\cJ_\cA$ or $\cA^h\cap X\in\cI_\cA\subseteq\cJ_\cA$.

Suppose that for some $\cA^h\in\cB[\cA]$ we have that $\cA^h\in\cJ_\cA$. Then $\cA^h\in\cJ_\cA\setminus\cI_\cA$ and the minimum $\alpha$ such that $\cA^h\in\cJ_\alpha$ is a successor ordinal. Let $\alpha=\beta+1$. By the definition of $\alpha$, we have that $\cJ_\alpha\neq\cJ_\beta$. Hence $\cJ_\alpha=\langle\cJ_\beta\cup\{\omega\setminus X_\beta\}\rangle$ and there is $u\in\cJ_\beta$ such that $\cA^h\subseteq u\cup(\omega\setminus X_\alpha)$. It follows that $\cJ_\alpha=\cJ_\beta$ by definition. A contradiction.

To see that it is complete, let $X\subseteq\omega$ such that $X$ is $\cJ_\cA$-big. Since $\cJ_\cA^+\subseteq\cI_\cA^+$, we see that $X$ is also $\cI_\cA$-big and there is $\alpha<\mathfrak{c}$ such that $X=X_\alpha$. Consider $\cJ_\alpha$. If there are $\cA^h\in\cB[\cA]$ and $u\in\cJ_\alpha$ such that
$\cA^h\subseteq u\cup (\omega\setminus X_\alpha)$, thence $\cA^h\cap X_\alpha\subseteq u\in\cJ_\alpha\subseteq\cJ_\cA$, contradicting that $X$ is $\cJ_\cA$-big. Henceforth, $\cJ_{\alpha+1}=\langle\cJ_\alpha\cup\{\omega\setminus X_\alpha\}\rangle$ and $(\omega\setminus X)\in\cJ_\cA$.
\end{proof}

\begin{teo}\label{maintheorem}
It is consistent that $\mathfrak{i}_a<\mathfrak{c}$.
\end{teo}

\begin{proof}
Let $V$ be a model of $\mathfrak{c}=\omega_2$.
We will recursively define a finite support iteration of forcing notions $\langle\mathbb{P}_\alpha,\dot{\mathbb{Q}}_\beta:\alpha\leq\omega_1,\beta<\omega_1\rangle$ and an increasing sequence of (names of) independent families $\{\dot{\cA}_\alpha:\alpha<\omega_1\}$, each accompanied with a complete diagonalization ideal $\dot{\cJ}_\alpha$.

$\underline{\alpha=0}$. Let $\cA_0$ be any countable independent family which separates points. Let $\cJ_0$ be a complete diagonalization ideal for $\cA_0$ and define $\mathbb{P}_0=\emptyset$.

$\underline{\alpha=\xi+1}$. Let $\dot{\mathbb{Q}}_\xi$ be a $\mathbb{P}_\xi$-name for the poset $\mathbb{M}(\cJ_\xi)$ and define $\mathbb{P}_\alpha=\mathbb{P}_\xi\ast\dot{\mathbb{Q}}_\xi$. In $V^{\mathbb{P}_\alpha}$, let $A_\xi$ be the $\mathbb{M}(\cJ_\xi)$-generic real over $V^{\mathbb{P}_\xi}$ and let $\cJ_\alpha$ be a complete diagonalization ideal for $\cA_\alpha=\cA_\xi\cup\{A_\xi\}$.
Let $\dot{\cA}_\alpha$ and $\dot{\cJ}_\alpha$ be $\mathbb{P}_\alpha$-names for $\cA_\alpha$ and $\cJ_\alpha$ respectively.

$\underline{\alpha\textnormal{ limit}}$. Let $\mathbb{P}_\alpha$ be the finite support iteration $\langle\mathbb{P}_\beta,\dot{\mathbb{Q}}_\gamma:\beta\leq\alpha,\gamma<\alpha\rangle$.
Working in $V^{\mathbb{P}_\alpha}$, let $\cA_\alpha=\bigcup_{\beta<\alpha}\cA_\beta$ and let $\cJ_\alpha\supseteq\cI_\alpha\supseteq\bigcup_{\beta<\alpha}\cJ_\beta$ be a complete diagonalization ideal for $\cA_\alpha$. 
Hence let $\dot{\cA}_\alpha$ and $\dot{\cJ}_\alpha$ be $\mathbb{P}_\alpha$-names for $\cA_\alpha$ and $\cJ_\alpha$ respectively.

Let $G$ be a $\mathbb{P}_{\omega_1}$-generic filter.
In $V[G]$, it is clear  from theorem \ref{independence} that $\cA=\bigcup_{\alpha<\omega_1}\cA_\alpha$ is an independent family. 
We can index 
$$\cA=\{A'_n:n\in\omega\}\cup\{A_\alpha:\alpha<\omega_1\},$$ where $\cA_0=\{A'_n:n\in\omega\}$ and $A_\alpha$ is the generic real added at the $\alpha$-stage of the iteration for every $\alpha<\omega_1$. 
Since $\cA_0\subseteq\cA$, we have that $\cA$ separates points. Suppose that $X\subseteq\omega$ is $\cA$-dense, $n\in X$ and let $p\in\mathbb{P}_{\omega_1}$ be a condition that forces these two statements. 
Since $\mathbb{P}_{\omega_1}$ is ccc, there exists an $\alpha<\omega_1$ such that $X\in V^{\mathbb{P}_\alpha}$ and $supp(p)\subseteq\alpha$. 
Hence 
$$V[G\rest\alpha]\models(1_{\mathbb{M}(\cJ_\alpha)}\forces``X\textnormal{ is }\dot{\cA}\textnormal{-dense''}),$$
where $p\in G\rest\alpha$. Otherwise, there should be a condition $q\leq p$ such that 
$$q\forces``X\textnormal{ is not }\cA\textnormal{-dense''}.$$ 
Thus, applying lemma \ref{A-dense}, we see that $X$ is $\cJ_\alpha$-big, and since $\cJ_\alpha$ is a complete diagonalization ideal, it follows that $(\omega\setminus X)\in\cJ_\alpha$. In consequence, the generic real $A_\alpha$ is forced to be almost contained in $X$. Consider any condition $r=(s,F)\in\mathbb{M}(\cI_\alpha)$ such that$n\in s$. Let $p'\leq_\alpha p$ forcing this and define $q=(p')\conc r\in\mathbb{P}_{\alpha+1}\subseteq\mathbb{P_{\omega_1}}$. Notice that 
$$q\forces ``n\in A_\alpha\subseteq^* X\textnormal{''}.$$

By density, $\cA$ is an anchored independent family of size $\omega_1$, and since $\mathbb{P}_{\omega_1}$ is ccc, it preserves cardinalities and therefore $V[G]\models(\omega_1=\mathfrak{i}_a<\mathfrak{c}=\omega_2)$.

\end{proof}

\begin{cor}
It is consistent with {\sf{ZFC}} that there is a disjointly tight irresolvable space of weight less that $\mathfrak{c}$. 
\end{cor}

\begin{proof}
Use the anchored independent family constructed above. By proposition \ref{submaximalequivalence}, $2^{\omega_1}$ contains a countable dense subspace $X$ which is submaximal (and then irresolvable). By theorem \cite{bella2020disjointly}, $X$ is disjointly tight. Thus, since $\mathfrak{c}=\omega_2$ in the model of theorem \ref{maintheorem}, $w(X)=\omega_1<\mathfrak{c}$.
\end{proof}

\section{Remarks and questions}

In \cite{moore2004parametrized}, parametrized versions of $\diamondsuit$ are introduced for several cardinal invariants of the continuum. It is also proved that the study of parametrized diamond principles is, in some way, dual to the study of cardinal invariants of the continuum, in the sense that the corresponding parametrized diamond principle $\diamondsuit(\eta)$ implies $\eta=\omega_1$, alike $\diamondsuit$ implies {\sf CH}. On the other hand, all this principles are independent with {\sf CH}.
These parametrized versions of diamond are considered for definable cardinal invariants, however, they even have an impact on non-definable cardinals.
An example of this is that $\diamondsuit(\mathfrak{b})$ implies $\mathfrak{a}=\omega_1$ and that $\diamondsuit(\mathfrak{r})$ implies $\mathfrak{u}=\omega_1$ \cite{moore2004parametrized}. 
With respect to irresolvable spaces, it is also known that $\diamondsuit(\mathfrak{r}_\mathbb{Q})$ implies $\mathfrak{i}=\omega_1$ (see \cite{moore2004parametrized}) where $\mathfrak{r}_\mathbb{Q}$ is the reaping number of dense subsets of the rational numbers. In \cite{cancino2014countable}, the cardinal $\mathfrak{r}_{\textsf{scat}}$ is defined as the reaping number of the quotient algebra $\cP(\mathbb{Q})$ modulo the ideal of scattered subsets of the rationals and it is proved that $\diamondsuit(\mathfrak{r}_{}\textsf{scat})$ implies that $\mathfrak{irr}=\omega_1$. We do not know if the same is holds for $\mathfrak{sbm}$:

\begin{que}
Does $\diamondsuit(\mathfrak{r}_{\textsf{scat}})$ imply that $\mathfrak{sbm}=\omega_1$? Or even $\mathfrak{i}_a=\omega_1$?
\end{que}

A more general question is if we can answer the question posed by Bella and Hru\v{s}\'ak using parametrized diamond principles.

\begin{que}
Does it follow from a parametrized diamond principle that there is a countable, irresolvable space, that is disjointly tight and has weight $\omega_1$?
\end{que}

There may always be countable, dense, submaximal subspaces of $2^\mathfrak{i}$, since we do not know the answer to the following question. 

\begin{que}
Is it consistent with {\sf ZFC} that $\mathfrak{i}<\mathfrak{i}_a$?
\end{que}

Although the relationship between $\mathfrak{i}$ and countable irresolvable spaces is very natural, it is not known if there could be a countable irresolvable space of $\pi$-weight less than $\mathfrak{i}$. Hence, it is still possible to have $\mathfrak{irr}=\mathfrak{i}$. The same  is true for $\mathfrak{sbm}$ and $\mathfrak{i}_a$. It seems that both problems could be solved with the same ideas. 

\begin{que}
Is it consistent to have $\mathfrak{sbm}<\mathfrak{i}_a$?
Are $\mathfrak{sbm}$ and $\mathfrak{i}$ provable comparable?
\end{que}

In case all inequalities are possible, it would be interesting to have a model in which $\mathfrak{i}$, $\mathfrak{i}_a$, $\mathfrak{irr}$ and $\mathfrak{sbm}$ are all different.

\bibliographystyle{plain}
\bibliography{resolvability.bib}

\end{document}